\documentclass[reqno, 10pt]{article}
\usepackage{amssymb, amsmath, amsthm, amsfonts, amscd, epsfig}
\usepackage{color}

\usepackage[english]{babel}
\usepackage{bbm}
\usepackage{graphicx, epsfig}
\usepackage{geometry,graphicx,pgfplots,caption,booktabs}
\usepackage{verbatim}
\usepackage{dsfont}
\usepackage{bm}

\usepackage[titletoc,toc,title]{appendix}
\usepackage[normalem]{ulem}
\newtheorem{theorem}{Theorem}[section]

\newtheorem{lemma}{Lemma}[section]
\newtheorem{proposition}[theorem]{Proposition}
\newtheorem{definition}{Definition}[section]

{\bf}{\it}

\newcommand{\Om}{\Omega}

\newcommand{\dd}{\,\mathrm{d}}
\newcommand{\ip}[2]{\left\langle #1,#2\right\rangle_{\partial\Om}}
\newcommand{\norm}[1]{\left\lVert #1\right\rVert}
\newcommand{\R}{\mathbb R}

\newcommand{\p}{\partial}

\newcommand{\beq}{\begin{equation}}
	\newcommand{\eeq}{\end{equation}}
\newcommand{\RN}[1]{%
	\textup{\uppercase\expandafter{\romannumeral#1}}%
}

\numberwithin{equation}{section}
\numberwithin{figure}{section}

\begin{document}
	
	\title{Recovering Small Perturbations of a Known Inclusion from One Mixed Steklov Eigenpair with Local Boundary Data}
	%Stable Reconstruction of Small Dirichlet-Inclusion Perturbations from One Steklov--Dirichlet Eigenpair
	
	\author{Jiho Hong\thanks{Department of Mathematics, The Chinese University of Hong Kong, Shatin, N. T., Hong Kong SAR, P.R. China (\texttt{jihohong@alumni.kaist.ac.kr})}}
	
	\maketitle
	
	\begin{abstract}
		We study the recovery of small normal displacements of a known inclusion from a single mixed Steklov eigenpair. The data consist of the eigenvalue and the eigenfunction trace on any nonempty open subset of the fixed outer boundary, with the eigenfunction normalized using only that subset. For an inner Dirichlet condition, we prove that the linearized measurement map is injective at every simple eigenvalue. A compatible adjoint identity yields a reconstruction method with a conditional stability estimate on any fixed finite-dimensional space of smooth displacements. For the principal Dirichlet eigenpair, the same local data also determine the inclusion uniquely among sufficiently small normal graphs. In a three-dimensional concentric shell, the eigenvalue determines the mean displacement at first order, while spherical harmonics diagonalize the remaining full-boundary trace map. The singular values decay exponentially with harmonic degree, yielding logarithmic conditional stability under a Sobolev prior. Numerical experiments using the perturbed geometry illustrate the quadratic linearization error and the deterioration of conditioning as the observed patch shrinks. For an inner Neumann condition, we prove linearized injectivity at every positive simple eigenvalue.
\end{abstract}

%\tableofcontents

\section{Introduction and Main Results}\label{sec:main}

Let $\Om$ and $D$ be nonempty open and bounded domains in $\R^d$ with $d\geq3$ satisfying $\overline D\subset\Om$. We assume that $A:=\Om\setminus\overline D$ is connected and both $\partial D$ and  $\partial\Om$, inner and outer boundaries, are $C^{3,\alpha}$ for some
$0<\alpha<1$.
Let $\mathsf{n}_D$ be the unit normal vector pointing out of $D$ and assume that $h\in C^{3,\alpha}(\partial D)$.
We define perturbed bounded domains $D_\varepsilon$ and $A_\varepsilon$ as
\begin{equation}\label{eq:graph}
	\partial D_\varepsilon=\{x+\varepsilon h(x)\mathsf{n}_D(x)\,:\,x\in\partial D\}\quad\mbox{and}\quad A_\varepsilon=\Om\setminus\overline{D_\varepsilon}\quad\forall \varepsilon\in (-\varepsilon_0,\varepsilon_0)
\end{equation}
for sufficiently small $\varepsilon_0>0$. We write \(p=\varepsilon h\)
for the physical normal displacement on \(\partial D\).
Let $\nu$ be the outward unit normal of $A_\varepsilon$, so that $\nu=\mathsf{n}_\Om$ on $\partial\Om$ and $\nu=-\mathsf{n}_{D_\varepsilon}$ on $\partial D_\varepsilon$.
We consider eigenpairs $(\lambda_\varepsilon,u_\varepsilon)\in\R\times (H^1(A_\varepsilon)\backslash\{0\})$ of the mixed Steklov--Dirichlet problem
\begin{equation}
	\label{eq:evp}
	\left\{
	\begin{aligned}
		-\Delta u_\varepsilon&=0\quad\mbox{in }A_\varepsilon,\\
		\partial_\nu u_\varepsilon&=\lambda_\varepsilon u_\varepsilon\quad\mbox{on }\partial\Omega,\\
		u_\varepsilon&=0\quad\mbox{on }\partial D_\varepsilon.
	\end{aligned}
	\right.
\end{equation}
The weak formulation and the compact trace embedding give a discrete
positive spectrum, where positivity follows from the energy identity $\|\nabla u_\varepsilon\|_{L^2(A_\varepsilon)}^2
=\lambda_\varepsilon \|u_\varepsilon\|_{L^2(\partial\Om)}^2$ and the Dirichlet condition on the nonempty inner boundary.
We aim to quantify what a single measured mode on an accessible boundary patch can reveal about the inner boundary that cannot be observed directly.
The Steklov condition converts the measured eigenpair into local Cauchy data.
The challenges are ill-posedness and nonlinearity of the inverse problem:
Boundary unique continuation determines the field from exact data \cite{Tolsa:2023:UCP}, but continuing noisy data toward the inclusion is ill-posed \cite{Alessandrini:2000:OS}; Perturbed inclusion changes both the eigenvalue and its trace nonlinearly.

In a quasistatic conduction interpretation, the harmonic field is the voltage, and the inner Dirichlet and Neumann conditions model a grounded conductor and an insulating cavity, respectively.
The outer Steklov condition selects voltage patterns whose normal derivative is proportional to the voltage.
The same Steklov--Neumann boundary coupling arises in linearized fluid sloshing, with the spectral condition on the free surface and zero normal flow at rigid walls \cite{BKPS2010,LPPS2022Sloshing}.
Studies of mixed Steklov eigenvalue bounds, shape optimization, and corner asymptotics show how the spectrum depends on geometry and boundary conditions \cite{GPPS2022Hole,HL2019Mixed,HLL:2025:ASD,HLS:2026:SDE,LPPS2022Corners}.
Steklov spectral geometry and its Dirichlet-to-Neumann interpretation
are surveyed in~\cite{GP,CGGS}. Variational spectral continuity and shape optimization in~\cite{GKL} concern general geometric limits; The shrinking-hole limits were considered in~\cite{CPN,Nazarov}.

We investigate the inverse problem to recover the perturbed Dirichlet inclusion $D_\varepsilon$ from the reference domain $A=A_0$ and one Steklov--Dirichlet eigenpair.
The eigenfunction trace carries spatial information beyond the scalar eigenvalue.
We prove that one simple eigenvalue and its eigenfunction trace on a nonempty open subset of $\partial\Omega$ determine the normal displacement from a known Dirichlet inclusion in both the original nonlinear problem (Theorem \ref{thm:nonlinear:unique}) and the linearized problem (Theorem \ref{theorem:uniqueness}).
For a Neumann inclusion, we prove the linearized identifiability at every positive simple eigenvalue (Theorem \ref{thm:injectivity}).
The linearized equations give a reconstruction method. We establish a conditional stability estimate for finite-dimensional spaces of smooth displacements (Proposition \ref{prop:local}).
For a concentric reference shell, we leverage the explicit expression of the reference principal eigenpair to prove that the linearized inverse problem is diagonalized (Proposition \ref{prop:diag}) and prove the infinite dimensional stability (Theorem \ref{thm:stability}).
Finally, we numerically test how the size of the local observation patch in $\partial\Omega$, the scale of the displacement and measurement noise affect the reconstruction.

The analysis on the linearization of the inverse problem has three main parts. First, under a fixed-domain formulation, we derive the shape derivatives and quadratic remainders for simple eigenpairs (Theorem \ref{thm:simple}).
Because we establish the local analyticity of the map from the extended pull-back transformation (rather than just $\varepsilon$) to the eigenpair (Lemma \ref{lem:analytic}), the constants in the remainder are uniform in both $\varepsilon$ and $h$.
Second, a compatible adjoint identity converts the perturbed boundary information of eigenfunction into displacement moments (Theorem \ref{thm:adjoint}). Along with the moments for perturbed eigenvalues, they give local injectivity and finite-dimensional conditional stability for the simple mode. Third, restricting to the case when $A$ is a concentric annulus, explicit spherical coefficients quantify the exponential loss of high-degree information and motivate regularization.
The local-data theorem and the full-boundary spherical stability theorem have different scopes and are stated separately.

Write $H_D^1(A)=\{v\in H^1(A):\gamma_D v=0\}$, where $\gamma_D$ is the continuous trace operator $H^1(A)\to H^{1/2}(\partial D)$. Because \(h\in C^{3,\alpha}(\partial D)\), for sufficiently small $|\varepsilon|$, one can choose $V\in W^{2,\infty}(\mathbb{R}^d,\mathbb{R}^d)$ equal to $h \mathsf{n}_D$ on $\partial D$, zero in a neighborhood of $\partial\Om$ and $(I+\varepsilon V)|_A$ is a diffeomorphism onto $A_\varepsilon$.
Put $F_\varepsilon:=I+\varepsilon V$.
The material pullback $U_\varepsilon:=u_\varepsilon\circ F_\varepsilon$ allows us to compare the eigenfunctions on the fixed domain $A$.
The following theorem gives first-order expansions of $\lambda_\varepsilon$ and $U_\varepsilon$ with quadratic remainders.

\begin{theorem}[Shape derivative with local boundary normalization]
	\label{thm:simple}
	Let \((\lambda_0,u_0)\) be a simple Steklov--Dirichlet eigenpair on \(A\),
	with \(\|u_0\|_{L^2(\partial\Omega)}=1\). Let \(\Gamma\) be any nonempty
	relatively open subset of \(\partial\Omega\), and put \(m_\Gamma=\|u_0\|_{L^2(\Gamma)}^2\). Write
	\(\gamma_\Omega\) for the outer trace and set
	\(\ell_\Gamma(v)=\int_\Gamma(\gamma_\Omega v)u_0\,d\sigma\) for
	\(v\in H^1(A)\). In particular, \(\ell_\Gamma(u_0)=m_\Gamma\).
	Then \(m_\Gamma>0\). The nearby eigenpair branch \((\lambda_\varepsilon, u_\varepsilon)\) is uniquely fixed near \((\lambda_0,u_0)\), in pulled-back coordinates, by \(\langle u_\varepsilon, u_0\rangle_{L^2(\Gamma)}=m_\Gamma\).
	Set \( U_\varepsilon= u_\varepsilon\circ F_\varepsilon\) and define
	$\lambda'[h]=\int_{\partial D} h(\partial_\nu u_0)^2 \dd\sigma$.
	The Eulerian derivative \(v_\Gamma=u'_\Gamma[h]\) is the unique \(H^2(A)\) solution of
	\begin{equation}\label{eq:local-shape}
		\left\{
		\begin{aligned}
			\Delta v_\Gamma&=0 \quad\mbox{in } A,\\
			(\partial_\nu-\lambda_0)v_\Gamma&=\lambda'[h]u_0 \quad\mbox{on }\partial\Omega,\\
			v_\Gamma&=h\partial_\nu u_0 \quad\mbox{on }\partial D,
		\end{aligned}
		\right.
		\quad\mbox{subject to}\quad \langle v_\Gamma, u_0\rangle_{L^2(\Gamma)}=0.
	\end{equation}
	There exist constants \(C,\varepsilon_0>0\) such that, for \(|\varepsilon|<\varepsilon_0\),
	\begin{align}
		|\lambda_\varepsilon-\lambda_0-\varepsilon\lambda'[h]|
		&\le C\varepsilon^2, \label{eq:eigen-remainder}\\
		\|  U_\varepsilon-u_0-\varepsilon(u'_\Gamma[h]+V\cdot\nabla u_0)\|_{H^1(A)}
		&\le C\varepsilon^2, \label{eq:material-remainder}\\
		\|(  u_\varepsilon-u_0-\varepsilon u'_\Gamma[h])|_\Gamma\|_{H^{1/2}(\Gamma)}
		&\le C\varepsilon^2. \label{eq:local-remainder}
	\end{align}
	For fixed $A$, $\Gamma$ and $(u_0,\lambda_0)$, the constants can be
	chosen uniformly over \(h\) if the extensions \(V\) have uniformly bounded \(W^{2,\infty}\) norms and vanish on a common neighborhood of $\partial\Omega$.
\end{theorem}
The \(C^{3,\alpha}\) boundary assumptions suffice for the quadratic remainder and the trace differentiations.

We define the compatible test space and the adjoint complement by
$$\mathcal G_\Gamma=\{g\in L^2(\Gamma)\,:\,\langle{g},{u_0}\rangle_{L^2(\Gamma)}=0\}\quad\mbox{and}\quad \mathcal{H}_\perp=\{ w\in H_D^1(A)\,:\,\langle w,u_0 \rangle_{L^2(\partial\Om)}=0\}.$$

\begin{theorem}[Compatible adjoint identity for local boundary data]
	\label{thm:adjoint}
	Let the assumptions of Theorem \ref{thm:simple} hold. Denote by \(\widetilde{g}\in L^2(\partial\Omega)\) the zero extension of \(g\in\mathcal G_\Gamma\).
	For all \(g\in\mathcal G_\Gamma\), there is a unique weak solution
	\(w_g\in \mathcal{H}_\perp\) of
	\begin{equation}\label{eq:adjoint}
		\left\{\begin{aligned}
			\Delta w_g&=0\quad \text{in }A,\\
			(\partial_\nu-\lambda_0)w_g&=\widetilde{g}\quad \text{on }\partial\Omega,\\
			w_g&=0\quad\text{on }\partial D.
		\end{aligned}\right.
	\end{equation}
	The weak solution satisfies the following estimate:
	\[
	\|w_g\|_{H^1(A)}+\|\partial_\nu w_g\|_{H^{1/2}(\partial D)}
	\le C\|g\|_{L^2(\Gamma)}.
	\]
	There exist constants \(C,\varepsilon_0>0\) such that the following identities and bounds hold for \(p=\varepsilon h\) and \(|\varepsilon|<\varepsilon_0\).
	\begin{align}
		\int_{\partial D} p(\partial_\nu u_0)(\partial_\nu w_g) \dd\sigma
		&=-\int_\Gamma(  u_\varepsilon-u_0)g \dd\sigma
		+R_\varepsilon(g),\qquad
		|R_\varepsilon(g)|\le C\varepsilon^2\|g\|_{L^2(\Gamma)}, \label{eq:patch-moment}\\
		\int_{\partial D} p(\partial_\nu u_0)^2 \dd\sigma
		&=\lambda_\varepsilon-\lambda_0+r_\varepsilon,\qquad
		|r_\varepsilon|\le C\varepsilon^2. \label{eq:eigen-moment}
	\end{align}
	The constants $C$ are independent of $\varepsilon$ and \(g\), with the uniformity in \(h\) of Theorem~\ref{thm:simple}.
\end{theorem}

The solution operator $g\mapsto w_g$ of the adjoint problem is bounded from $\mathcal G_\Gamma$ to $\mathcal H_\perp$.
In a fixed neighborhood of $\partial D$, separated from $\partial\Omega$, the harmonic function $w_g$ is $H^2$ by homogeneous Dirichlet regularity and $\partial_\nu w_g|_{\partial D}\in H^{1/2}(\partial D)\subset L^2(\partial D)$,
with norm controlled by $\norm g_{L^2(\Gamma)}$.
Also, \(u_0\) belongs to \(H^3(A)\cap C^{2,\alpha}(\overline A)\) by the regularity argument in Section~\ref{sec:spaces}.
Consequently, $(\partial_\nu u_0)(\partial_\nu w_g)\in L^2(\partial D)$, and the moments on the left-hand sides of \eqref{eq:patch-moment} and \eqref{eq:eigen-moment} are well-defined for $p\in L^2(\partial D)$.
The next theorem establishes linearized identifiability from the local data of one simple eigenpair.
\begin{theorem}[Linearized Steklov--Dirichlet identifiability]
	\label{theorem:uniqueness}
	Let the assumptions of Theorem~\ref{thm:simple} hold.
	For all \(g\in\mathcal G_\Gamma\), let \(w_g\in \mathcal{H}_\perp\) be the unique weak solution of \eqref{eq:adjoint}.
	Suppose that $h\in C^{3,\alpha}(\partial D)$ satisfies
	\begin{equation}\label{eq:h:system}
		\int_{\partial D} h\, (\partial_\nu u_0)^2\dd\sigma=0\quad\mbox{and}\quad\int_{\partial D} h(\partial_\nu u_0) (\partial_\nu w_g)\dd\sigma=0\quad\forall g\in\mathcal{G}_\Gamma.
	\end{equation}
	Then $h$ is identically zero.
\end{theorem}
The linear moments in \eqref{eq:h:system} can be estimated from the
measurements via \eqref{eq:patch-moment} and \eqref{eq:eigen-moment},
with quadratic remainders. Proposition~\ref{prop:local} uses these
moments to obtain a finite-dimensional stability estimate.

Theorem~\ref{thm:nonlinear:unique} gives nonlinear uniqueness for the principal Dirichlet mode.
Theorem~\ref{thm:injectivity} establishes linearized identifiability for positive simple Steklov--Neumann eigenpairs.
Theorem~\ref{thm:stability} quantifies the resolution of full-boundary reconstruction in a three-dimensional concentric shell.

\paragraph{Related work and organization.}
Classical multidimensional Borg--Levinson theory recovers an interior potential from Dirichlet eigenvalues and boundary normal derivatives of the associated eigenfunctions~\cite{NSU}. 
On mixed Steklov eigenpairs, earlier work derived shape derivatives for translating a Dirichlet hole \cite{HLS,GP:2024:mono}, analyzed mixed Steklov spectra as holes shrink \cite{CPN}, and extracted mixed Steklov eigenvalues from scattering data as signatures of thin screens \cite{Cakoni:2022:TSTS}. The closest methodological precedent uses modal perturbation formulas and a dual identity to recover changes in a conductivity transmission interface~\cite{ABFKL}. 
Prior inverse Steklov work has also recovered polygonal geometry from high-frequency spectral asymptotics \cite{Krymski:2021:ISSP}.

Besides our new inverse problem setting, our new step of methodology is to use compatible tests supported on the observation patch to solve an adjoint problem at the reference eigenvalue and convert the measured trace change into moments of the inner-boundary displacement. Together with the eigenvalue shift, these moments give a least-squares reconstruction with a quadratic linearization remainder and a conditional stability estimate on fixed finite-dimensional displacement spaces. We prove linearized injectivity for every simple Dirichlet mode and every positive simple Neumann mode, and nonlinear uniqueness for the principal Dirichlet mode among sufficiently small normal graphs. For the linearized full-boundary problem in a three-dimensional concentric shell, explicit singular values quantify exponential attenuation of the higher harmonics and yield logarithmic conditional stability under a Sobolev prior. These inverse conclusions complement the forward hole-shape optimization of Fern\'andez Bonder et al.~\cite{BGR}
and the shrinking-hole spectral asymptotics of~\cite{CPN,Nazarov}.

The analytic framework uses spectral perturbation~\cite{Kato},
boundary perturbations~\cite{Henry}, and domain calculus~\cite{HenrotPierre,DelfourZolesio}. The trace and potential spaces follow~\cite{McLean,AK,AKL}.
%; contour formulas use~\cite{GS}.
Kirsch's obstacle-scattering shape derivative~\cite{Kirsch} provides a methodological comparison with a different observation operator. For stability, elliptic Cauchy estimates~\cite{ARRV} and inverse-inclusion instability~\cite{DCR} motivate the distinction between uniqueness and recoverability. Our explicit rate follows from the calculated spherical coefficients, within the regularization frameworks of~\cite{EHN,Hansen}; those other inverse problems do not supply a stability theorem for the present unknown boundary.

Section~\ref{sec:proof} develops the functional framework and spectral
derivatives. Section~\ref{sec:inverse} gives the adjoint reconstruction
and local-data result.
Section~\ref{sec:neumann} discusses the inner Neumann extension.
Section~\ref{sec:sphere} treats spherical coefficients and stability, and Section~\ref{sec:numerics} presents the numerical tests.
Section~\ref{sec:scope} states the conclusions.

\section{Shape Differentiation of Simple Steklov–Dirichlet Eigenpairs}\label{sec:proof}
\subsection{Energy and Trace Spaces}\label{sec:spaces}

We use standard Sobolev trace and conormal derivative constructions~\cite{McLean}.
Unless complex spectral parameters are used, all function spaces are real in this paper.
The spaces $H^s(G)$ are the $L^2$-based Sobolev spaces on a domain or
compact boundary $G$; on a boundary they are defined using local charts.
On a closed boundary $\Sigma$, we write
$H^{-s}(\Sigma)=(H^s(\Sigma))^*$ with $L^2(\Sigma)$ as pivot, for $s>0$.
Negative Sobolev exponents below are used on boundaries.

For $v\in H^1(A)$, let $\gamma_\Om v$ and $\gamma_Dv$ denote its
traces on the two boundary components. The trace maps are bounded into
$H^{1/2}(\partial\Om)$ and $H^{1/2}(\partial D)$, respectively. The energy space is
\[H_D^1(A):=\{v\in H^1(A)\,:\,\gamma_D v=0\}\]
equipped with the usual $H^1(A)$ norm.
By Poincare's inequality, the norm is equivalent to the gradient norm.
The corresponding energy space for perturbed domain is $H_D^1(A_\varepsilon)$.
For $u\in H^1(A)$ harmonic in distributions, the weak conormal derivative
$\partial_\nu u\in H^{-1/2}(\partial A)$ is defined by
\begin{equation}\label{eq:weaknormal}
	\langle\partial_\nu u,\gamma v\rangle_{H^{-1/2},H^{1/2}}
	=\int_A\nabla u\cdot\nabla v\dd x,\qquad v\in H^1(A).
\end{equation}
It is well-defined because the right side vanishes on $H_0^1(A)$.
Since $\partial\Om$ and $\partial D$ are separated, we have $\partial_\nu u|_{\partial\Om}\in H^{-1/2}(\partial \Om)$ and $\partial_\nu u|_{\partial D}\in H^{-1/2}(\partial D)$.

Thus the weak eigenproblem means $u_\varepsilon\in H_D^1(A_\varepsilon)\setminus\{0\}$ and
\begin{equation}\label{eq:weak:problem}
	\int_{A_\varepsilon}\nabla u_\varepsilon\cdot\nabla v\dd x
	=\lambda_\varepsilon\int_{\partial\Om}
	(\gamma_\Om u_\varepsilon)(\gamma_\Om v)\dd\sigma,
	\qquad\forall v\in H_D^1(A_\varepsilon).
\end{equation}
The outer Robin equation holds in $H^{-1/2}(\partial\Om)$, with
$\gamma_\Om u_\varepsilon\in H^{1/2}(\partial\Om)$ identified with an
element of that dual space through the $L^2$ pairing. The weak problem requires only Lipschitz regularity for $\partial A_\varepsilon$. For the shape derivatives and quadratic remainders we introduce the stronger $C^{3,\alpha}$ hypotheses.

An extension $V\in
W^{2,\infty}(\R^d;\R^d)$, supported away from $\partial\Om$ and with
the specified inner trace, suffices for the pullbacks used here. In
particular $F_\varepsilon$ and its inverse act boundedly on $H^1$.
The \(H^2\) boundary estimates give \(u_0\in H^2(A)\), since the inner
Dirichlet data vanish and the outer Neumann data \(\lambda_0u_0\)
belong to \(H^{1/2}(\partial\Omega)\). The improved trace gives
\(\lambda_0u_0\in H^{3/2}(\partial\Omega)\), so the next-order
estimates yield \(u_0\in H^3(A)\). Boundary Schauder estimates
for the homogeneous Dirichlet and Robin conditions in separate
collars, together with interior harmonic regularity, give
\(u_0\in C^{2,\alpha}(\overline A)\). Hence
$\partial_\nu u_0|_{\partial D}\in H^{3/2}(\partial D)\cap C^{1,\alpha}(\partial D)$.
For smooth $h$, the derivative $u'[h]\in H^2(A)$ has inner trace
$h \partial_\nu u_0|_{\partial D}\in H^{3/2}(\partial D)$ and outer conormal derivative in
$H^{1/2}(\partial\Om)$. Its material counterpart
$\dot U=u'[h]+V\cdot\nabla u_0$ belongs to $H_D^1(A)$, whereas
$u'[h]$ generally does not have zero inner trace.
The construction of material and Eulerian derivatives in these spaces
is consistent with the domain-calculus frameworks
of~\cite{Henry,HenrotPierre,DelfourZolesio}.

\begin{comment}
	The geometric displacement class is $C^{3,\alpha}(\partial D)$ with small norm. Extending individual moment functionals continuously to $L^2(\partial D)$ does not assert that every $L^2$ function defines a boundary perturbation or an $H^1$ shape derivative. In particular, for $p\in L^2$ alone the trace $p \partial_\nu u_0|_{\partial D}$ need not be in $H^{1/2}$.
	For local data, restriction sends $H^{1/2}(\partial\Om)$ into $H^{1/2}(\Gamma)$, understood as the restriction space with its quotient norm. Our compactly supported $L^2$ tests avoid any ambiguity about extension by zero at the edge of $\Gamma$.    
\end{comment}

We use the following boundary continuation lemma in the identifiability proofs.
\begin{lemma}[Unique continuation, {\cite{AE:1997:UCP} and \cite[Theorem 2.3]{KN:1998:UCP}}]
	\label{lemma:UCP}
	Let $A'$ be a connected and bounded open subset of $\mathbb{R}^d$ with $C^{1,\alpha}$ boundary for some $\alpha\in(0,1]$.
	Let $\Gamma'$ be a nonempty and relatively open subset of $\partial A'$.
	If \(w\) is weakly harmonic in \(A'\), belongs to \(H^1\) locally
	up to \(\Gamma'\), and has zero Dirichlet trace and weak conormal
	trace there in the sense of \eqref{eq:weaknormal}, then \(w=0\) in \(A'\).
\end{lemma}
Tolsa~\cite[Theorem 1.1]{Tolsa:2023:UCP} treats the case of \(C^1\) boundaries. Note that our $C^{3,\alpha}$ assumption for $\partial A$ and $h$ is sufficient for use of Lemma \ref{lemma:UCP} and further used in establishing spectral differentiability and controlling remainders as described in the previous section.

\subsection{Shape Derivatives of Simple Eigenpairs}

We first establish analyticity of the pulled-back form, then use it to prove Theorem~\ref{thm:simple}.
\begin{definition}[{\cite[Definition 15.1]{Deimling:1985:NFA}}]
	Let $X$ and $Y$ be Banach spaces and let $\mathcal{U}$ be a nonempty open subset of $X$.
	The map \(\mathcal A:\mathcal{U}\to Y\) is real analytic if for every \(W_*\in\mathcal U\) there exist \(r>0\) and continuous symmetric \(n\)-linear maps \(B_n:X^n\to Y\), \(n\in\mathbb{N}\), such that $\sum_{n=1}^\infty\|B_n\|r^n<\infty$ and
	\[
	\mathcal A(W_*+K)=\mathcal A(W_*)+
	\sum_{n=1}^{\infty}B_n(K,\ldots,K),
	\qquad \|K\|_X<r.
	\]
\end{definition}
Let \(X\) be the closed real subspace of
\(W^{2,\infty}(\mathbb R^d;\mathbb R^d)\) consisting of vector fields
that vanish in a fixed collar of \(\partial\Omega\). For \(W\) in a sufficiently small open ball \(\mathcal U\subset X\) about zero, put
\[M(W)=I+DW,\qquad
C(W)=\det M(W)\,M(W)^{-1}M(W)^{-T}.\]
Define bounded symmetric linear forms $a_W,b:H_D^1(A)\times H_D^1(A)\to\R$ as
\[a_W(z,v)=\int_A C(W)\nabla z\cdot\nabla v\dd x\quad\mbox{and}\quad b(z,v)=\int_{\partial\Om}(\gamma_\Om z)(\gamma_\Om v)\dd\sigma\quad\forall z,v\in H_D^1(A).\]
Let $(\lambda(W),u(W))$ be a Steklov--Dirichlet eigenpair for Dirichlet inclusion being $(I+W)D$ so that $u(0)=u_0$.
Let $U(W)=u(W)\circ(I+W)$.
Via the change of variables with $I+W$, the variational identity \eqref{eq:weak:problem} for the Dirichlet inclusion $(I+W)D$ is equivalent to the pulled-back equation
\begin{equation}\label{eq:forms}
	a_W(U(W),v)
	=\lambda(W)b(U(W),v)\quad\forall v\in H_D^1(A).
\end{equation}
Define \(\mathcal A(W)\in\mathcal L(H_D^1(A),(H_D^1(A))^*)\) by \((\mathcal A(W)z)(v)=a_W(z,v)\) for all $z,v\in H_D^1(A)$.

\begin{lemma}
	The map \(\mathcal A:\mathcal{U}\to L(H_D^1(A),(H_D^1(A))^*)\) is real analytic.
\end{lemma}
\begin{proof}
	Equip \(L^\infty(A;\mathbb R^{d\times d})\) with the essential supremum of the Euclidean matrix operator norm. This is a Banach algebra under pointwise matrix multiplication. Choose \(\mathcal U\) so that \(\|DW\|_{L^\infty}<1\) throughout \(\mathcal U\).
	Fix \(W_*\in\mathcal U\) and set \(M_*=I+DW_*\).
	Then \(M_*^{-1}\in L^\infty(A;\mathbb R^{d\times d})\) by the Neumann series. For \(K\) small enough that \(W_*+K\in\mathcal U\) and
	\(\|M_*^{-1}DK\|_{L^\infty}<1\), the factorization
	\[
	M_*+DK=M_*(I+M_*^{-1}DK)
	\]
	gives
	\[
	(M_*+DK)^{-1}
	=\sum_{j=0}^\infty(-M_*^{-1}DK)^j M_*^{-1}.
	\]
	The \(j\)th term is a continuous homogeneous polynomial of degree
	\(j\) in \(K\), since \(K\mapsto DK\) is bounded linear and matrix
	multiplication is continuous. Its \(L^\infty\) norm is at most
	\[
	\|M_*^{-1}\|_{L^\infty}
	\bigl(\|M_*^{-1}\|_{L^\infty}\|DK\|_{L^\infty}\bigr)^j.
	\]
	This geometric majorant yields absolute convergence in \(L^\infty\)
	on a sufficiently small \(X\)-ball. Transposition is bounded linear,
	so the transposed inverse has the same property.
	The determinant \(\det(M_*+DK)\) is a finite sum of continuous
	homogeneous polynomials in \(K\) of degrees at most \(d\).
	The Cauchy product of these expansions converges absolutely in
	\(L^\infty\) by the submultiplicative norm bound. It gives
	\[
	C(W_*+K)=C(W_*)+\sum_{n=1}^\infty P_n(K),
	\]
	where \(P_n\) is a continuous homogeneous polynomial of degree \(n\)
	obtained from bounded multilinear products. Symmetrizing the
	multilinear map underlying \(P_n\) preserves its value on the
	diagonal \((K,\ldots,K)\).
	
	For a matrix coefficient \(C\in L^\infty(A;\mathbb R^{d\times d})\),
	define \(T_C:H_D^1(A)\to (H_D^1(A))^*\) by
	\[ (T_C z)(v)=\int_A C\nabla z\cdot\nabla v\dd x\quad\forall z,v\in H_D^1(A).\]
	Cauchy--Schwarz gives
	\[
	|(T_Cz)(v)|
	\le\|C\|_{L^\infty}\|\nabla z\|_{L^2(A)}\|\nabla v\|_{L^2(A)}
	\le\|C\|_{L^\infty}\|z\|_{H_D^1(A)}\|v\|_{H_D^1(A)}\quad\forall z,v\in H_D^1(A).
	\]
	Consequently \(C\mapsto T_C\) is bounded linear and
	\(\|T_C\|_{\operatorname{op}}\le\|C\|_{L^\infty}\).
	Applying it term by term gives
	\[
	\mathcal A(W_*+K)=\mathcal A(W_*)+
	\sum_{n=1}^{\infty}T_{P_n(K)},\qquad
	\sum_{n=1}^{\infty}\|T_{P_n(K)}\|_{\operatorname{op}}
	\le\sum_{n=1}^{\infty}\|P_n(K)\|_{L^\infty}<\infty .
	\]
	The symmetric multilinear coefficients of \(P_n\), composed with
	\(C\mapsto T_C\), supply the maps \(B_n\) in the definition.
	Since \(W_*\) was arbitrary, \(\mathcal A\) is real analytic in
	operator norm throughout \(\mathcal U\).
\end{proof}

\begin{lemma}\label{lem:analytic}
	For $W$ in an open neighborhood of zero in $X$, the pulled-back eigenpair branch $(U(W),\lambda(W))\in H_D^1(A)\times\mathbb R$ is uniquely fixed and real analytic as a function of $W$.
\end{lemma}
\begin{proof}
	The residual to which the theorem applies is
	\[
	\mathcal F(W,U,\lambda)
	=\bigl(a_W(U,\cdot)-\lambda b(U,\cdot),\
	\ell_\Gamma(U)-m_\Gamma\bigr)
	\in (H_D^1(A))^*\times\mathbb R .
	\]
	It is real analytic on an open neighborhood of \((0,u_0,\lambda_0)\)
	in \(X\times H_D^1(A)\times\mathbb R\), and
	\(\mathcal F(0,u_0,\lambda_0)=0\).

	Below, we prove that the derivative in \((U,\lambda)\) is a bounded isomorphism.
	Write \(L=a_0-\lambda_0 b:H_D^1(A)\to (H_D^1(A))^*\). Differentiation in the unknowns gives
	\begin{equation}\label{eq:DFzeta}
		D_{(U,\lambda)}\mathcal F(0,u_0,\lambda_0)[z,\eta]
		=\bigl(Lz-\eta b(u_0,\cdot),\ \ell_\Gamma(z)\bigr)\quad\forall (z,\eta)\in H_D^1(A)\times\mathbb{R}.
	\end{equation}
	
	First, we identify the range of \(L\). Poincar\'e's inequality on \(H_D^1(A)\) holds because \(A\) is connected and the Dirichlet boundary has
	positive surface measure. Thus \(a_0\) is coercive on \(H_D^1(A)\) and defines an isomorphism \(H_D^1(A)\to (H_D^1(A))^*\) by Lax--Milgram. The operator defined by
	\(b\) is compact, since it factors as \(T^*T\) with
	\(T:H_D^1(A)\to L^2(\partial\Omega)\) the compact outer trace. Therefore \(L\)
	is Fredholm of index zero. By simplicity \(\ker L=\mathbb Ru_0\), so
	its range is closed and has codimension one. Symmetry gives
	\((Lz)(u_0)=0\) for every \(z\in H_D^1(A)\). Since the functionals annihilating
	\(u_0\) also form a codimension-one subspace,
	\begin{equation}\label{eq:range:L}
		\operatorname{ran}L=\{f\in (H_D^1(A))^*\,:\, f(u_0)=0\}.
	\end{equation}

	Now prescribe an arbitrary target \((f,s)\in (H_D^1(A))^*\times\mathbb R\) as the output of \eqref{eq:DFzeta}.
	\begin{enumerate}
		\item Evaluating \(Lz-\eta b(u_0,\cdot)=f\) at \(u_0\) gives
		\(-\eta=f(u_0)\), since \(b(u_0,u_0)=1\). This determines \(\eta\) uniquely.
		\item The functional \(f+\eta b(u_0,\cdot)\) annihilates \(u_0\).
		The range characterization \eqref{eq:range:L} therefore gives the existence of a solution \(z_*\) of \(Lz_*=f+\eta b(u_0,\cdot)\).
		\item Every solution of this first equation of \eqref{eq:DFzeta} has the form \(z_*+cu_0\). The second equation of \eqref{eq:DFzeta}, \(\ell_\Gamma(z)=s\), determines
		\[c=\frac{s-\ell_\Gamma(z_*)}{m_\Gamma}.\]
		The patch nonvanishing argument at the start of the proof of Theorem~\ref{thm:simple} gives \(m_\Gamma>0\) independently of the branch construction.
	\end{enumerate}
	
	The derivative is thus a bounded bijection between Banach spaces.
	Its inverse is bounded by the bounded inverse theorem.
	We note that there are two distinct pairings in the calculation. The {full-boundary} form \(b\) comes from the Steklov equation and determines \(\eta\). The {patch} functional \(\ell_\Gamma\) fixes the eigenfunction amplitude.

	The analytic implicit function theorem \cite[Theorem 15.3(a)]{Deimling:1985:NFA} now gives neighborhoods of
	\(0\) and \((u_0,\lambda_0)\) and a unique real-analytic map
	\[W\longmapsto (U(W),\lambda(W))\in H_D^1(A)\times\mathbb R\]
	in those neighborhoods satisfying \(\mathcal F(W,U(W),\lambda(W))=0\). Uniqueness here is local among
	nearby solutions of the weak equation with the stated patch pairing.
\end{proof}

\begin{proof}[Proof of Theorem \ref{thm:simple}]
	First, we prove $m_\Gamma>0$ by contradiction:
	If \(m_\Gamma=0\), the trace of \(u_0\) vanishes on \(\Gamma\), and the
	Steklov condition gives \(\partial_\nu u_0=0\) there. Applying the unique continuation property, Lemma \ref{lemma:UCP}, for the harmonic function $u_0$ defined in $A$, we obtain \(u_0=0\) in $A$. This contradicts that $u_0$ is an eigenfunction.

	We apply Lemma \ref{lem:analytic} along \(W=\varepsilon V\). Write
	\(\dot U=\left.\frac{d}{d\varepsilon}  U_\varepsilon\right|_0\)
	and \(\dot\lambda=\left.\frac{d}{d\varepsilon}\lambda_\varepsilon\right|_0\).
	Taylor's formula gives
	\begin{equation}\label{eq:taylor}
		\|  U_\varepsilon-u_0-\varepsilon\dot U\|_{H^1(A)}
		+|\lambda_\varepsilon-\lambda_0-\varepsilon\dot\lambda|
		\le C\varepsilon^2 .
	\end{equation}
	Indeed, the branch has a bounded second derivative on a sufficiently
	small ball of vector fields by continuity at zero. For \(V\) in a bounded
	set, the lines \(W=\varepsilon V\) lie in that ball for a common
	\(\varepsilon_0\), and their second derivatives are bounded by a constant
	times \(\|V\|_{W^{2,\infty}}^2\). A fixed inner collar and cutoff in
	tubular coordinates give uniformly bounded extensions for bounded \(h\),
	with a common outer collar on which they vanish.
	
	\medskip
	\noindent\textbf{The Eulerian equation and boundary data.}
	Set \(v_\Gamma=\dot U-V\cdot\nabla u_0\in H^1(A)\), using the \(H^2\)
	regularity of \(u_0\). On every compact subset of \(A\),
	\(F_\varepsilon^{-1}=I-\varepsilon V+O(\varepsilon^2)\) uniformly.
	The expansion of
	\(  u_\varepsilon=  U_\varepsilon\circ F_\varepsilon^{-1}\)
	therefore gives
	\((  u_\varepsilon-u_0)/\varepsilon\to v_\Gamma\) in local \(L^2\).
	In detail, the remainder in \eqref{eq:taylor} stays \(O(\varepsilon^2)\)
	in local \(L^2\) by the change of variables and the bounded Jacobian.
	Composition of \(\dot U\) with \(F_\varepsilon^{-1}\) converges locally in
	\(L^2\) to \(\dot U\), by smooth approximation and the same Jacobian bound.
	The difference quotient of \(u_0\circ F_\varepsilon^{-1}\) converges locally
	uniformly to \(-V\cdot\nabla u_0\) by interior smoothness of \(u_0\).
	For every \(\phi\in C_c^\infty(A)\), its support lies in \(A_\varepsilon\)
	for small \(|\varepsilon|\). Testing the difference quotient against
	\(\Delta\phi\) gives zero, so \(\Delta v_\Gamma=0\) in distributions.
	
	In the outer collar, \(v_\Gamma=\dot U\). For tests \(\phi\in H_D^1(A)\)
	supported in that collar, differentiating the pulled-back equation gives
	\[
	\int_A \nabla v_\Gamma\cdot\nabla\phi\,dx
	=\lambda_0 b(v_\Gamma,\phi)+\dot\lambda b(u_0,\phi).
	\]
	Such tests can realize every outer \(H^{1/2}\) trace by lifting and
	cutting off. The weak conormal definition then gives
	\((\partial_\nu-\lambda_0)v_\Gamma=\dot\lambda u_0\) in
	\(H^{-1/2}(\partial\Omega)\). Differentiating the patch pairing gives
	\(\ell_\Gamma(v_\Gamma)=\ell_\Gamma(\dot U)=0\).
	On \(\partial D\) we have \( \dot U=0\), because the entire pulled-back
	branch lies in \(H_D^1(A)\). The tangential gradient of \(u_0\) vanishes there,
	so \(V\cdot\nabla u_0=-h\partial_\nu u_0\) and \(  v_\Gamma=h\partial_\nu u_0\).
	
	For clarity, the standard boundary regularity used here applies locally
	in separate collars of the two boundaries. For \(u_0\), homogeneous inner
	Dirichlet data and outer Neumann data \(\lambda_0  u_0\in H^{1/2}\)
	give \(H^2\) regularity. Applying the estimates at the next order gives
	\(u_0\in H^3(A)\), hence \(\partial_\nu u_0\in H^{3/2}(\partial D)\). Local Dirichlet
	Schauder regularity at \(\partial D\) also gives bounded \(\partial_\nu u_0\).
	For \(v_\Gamma\), the inner datum \(h\partial_\nu u_0\) belongs to \(H^{3/2}\) by
	smooth multiplication. Its outer Neumann datum
	\(\lambda_0  v_\Gamma+\dot\lambda  u_0\) belongs to \(H^{1/2}\).
	The \(H^2\) Dirichlet and Neumann estimates in the respective collars,
	together with interior regularity, give \(v_\Gamma\in H^2(A)\).
	
	\medskip
	\noindent\textbf{Eigenvalue, uniqueness and error.}
	Green's second identity gives
	\[
	0=\int_{\partial A}(u_0\partial_\nu v_\Gamma-v_\Gamma\partial_\nu u_0) \dd\sigma
	=\dot\lambda\int_{\partial\Omega}u_0^2 \dd\sigma
	-\int_{\partial D}h(\partial_\nu u_0)^2 \dd\sigma.
	\]
	The \(\lambda_0\) terms on the outer boundary cancel, and the first term
	of the inner integrand vanishes because \(u_0=0\) there. The full-boundary
	norm of the reference \(u_0\) is one, so this gives
	\(\dot\lambda=\lambda'[h]\). The difference of two solutions of
	\eqref{eq:local-shape} belongs to \(\ker L\) and has zero patch pairing.
	Since \(\ell_\Gamma(u_0)=m_\Gamma>0\), the difference is zero.
	Substitution of \(\dot U=v_\Gamma+V\cdot\nabla u_0\) into
	\eqref{eq:taylor} proves \eqref{eq:eigen-remainder} and
	\eqref{eq:material-remainder}. Taking the outer trace, using
	\(F_\varepsilon=I\) and \(V=0\) near \(\partial\Omega\), and restricting
	to \(\Gamma\) proves \eqref{eq:local-remainder}.
\end{proof}

\section{Identifiability and Reconstruction from Local Boundary Data}\label{sec:inverse}
The full measured perturbation belongs to
$\R\times\{v\in H^{1/2}(\partial\Om):\ip{v}{u_0}=0\}$
in the global gauge. The reconstruction error model uses the weaker
$\R\times L^2(\partial\Om)$ norm and, after applying finitely many
bounded tests, a Euclidean data norm. It does not posit pointwise
measurements of a weak conormal derivative. On $\Gamma$, the trace
belongs to the restriction space $H^{1/2}(\Gamma)$; the restricted
conormal derivative acts on $\widetilde H^{1/2}(\Gamma)$, the closure
of $C_c^\infty(\Gamma)$ in the full-boundary $H^{1/2}$ norm.
This specifies the distributional meaning of local Cauchy data.
The eigenvalue expansion is retained as a separate moment in the reconstruction. For the concentric reference geometry it recovers the mean displacement, which the compatible adjoints cannot detect; see Proposition~\ref{prop:diag}.

\subsection{Linearized Identifiability from Local Data}
The normalization \(\langle u_\varepsilon,u_0\rangle_{L^2(\Gamma)}=m_\Gamma\) uses only the trace on \(\Gamma\).
For local Cauchy data, choose \(\phi\in C_c^\infty(\Gamma)\)
with \(\int_\Gamma u_\varepsilon\phi\,d\sigma\ne0\). The weak
Steklov condition gives
\[
\lambda_\varepsilon=
\frac{\langle\partial_\nu u_\varepsilon,\phi\rangle_\Gamma}
{\int_\Gamma u_\varepsilon\phi\,d\sigma}.
\]
The noise model includes error in the estimated eigenvalue.

\begin{proof}[Proof of Theorem \ref{thm:adjoint}]
	We adopt \(b\) and \(L\) from the proof of Theorem~\ref{thm:simple}. The weak equation is
	\(Lw_g=f_g\), where \(f_g(\phi)=\int_\Gamma g\, \phi \dd\sigma\).
	Cauchy--Schwarz and the trace bound give
	\(\|f_g\|_{(H_D^1(A))^*}\le C\|g\|_{L^2(\Gamma)}\). The condition \(f_g(u_0)=0\)
	is necessary by symmetry and sufficient by the range description of \(L\).
	All solutions differ by a multiple of \(u_0\). The condition \(b(w_g,u_0)=0\)
	fixes that multiple because \(b(u_0,u_0)=1\).
	The restriction of \(L\) to the closed subspace
	\(\{w\in H_D^1(A)\,:\, b(w,u_0)=0\}\) is a bounded bijection onto its closed range.
	Its bounded inverse yields \(\|w_g\|_{H^1(A)}\le C\|g\|_{L^2(\Gamma)}\).
	Tests compactly supported in \(A\) give harmonicity. For all tests in \(H_D^1(A)\),
	the conormal definition and surjectivity of the outer trace on \(H_D^1(A)\)
	then give the outer equation in \eqref{eq:adjoint}.
	
	In a fixed smaller collar of \(\partial D\), local homogeneous Dirichlet
	regularity bounds the \(H^2\) norm of \(w_g\) by its \(H^1(A)\) norm.
	The normal trace bound there proves the estimate on \(\partial_\nu w_g\).
	In particular \((\partial_\nu u_0)(\partial_\nu w_g)\in L^2(\partial D)\) since \(\partial_\nu u_0\) is bounded,
	so the moment extends continuously to \(L^2\) displacements.
	
	Apply Green's second identity to \(w_g\) and \(u_\Gamma'[h]\). By harmonicity, we have
	$$\langle\partial_\nu u_\Gamma'[h],w_g\rangle_{\partial A}
	=\int_A\nabla u_\Gamma'[h]\cdot\nabla w_g\dd x=\langle\partial_\nu w_g,u_\Gamma'[h]\rangle_{\partial A}.$$
	With the conormal paired against the trace on each boundary component, it gives
	\begin{align*}
		0&=\langle\partial_\nu u_\Gamma'[h],w_g\rangle_{\partial A}
		-\langle\partial_\nu w_g,u_\Gamma'[h]\rangle_{\partial A}\\
		&=\lambda'[h]\int_{\partial\Omega}w_g u_0 \dd\sigma
		-\int_\Gamma u_\Gamma'[h] g \dd\sigma
		-\int_{\partial D} h(\partial_\nu u_0)(\partial_\nu w_g) \dd\sigma .
	\end{align*}
	In the second line, the two outer \(\lambda_0\) terms cancel, the zero extension restricts the source term to \(\Gamma\), and on
	\(\partial D\) the trace of \(w_g\) vanishes while \(u_\Gamma'[h]|_{\partial D}=h\partial_\nu u_0\).
	The full-boundary adjoint constraint kills the term proportional to
	\(\lambda'[h]\). This pairing uses the reference solution and the computed
	adjoint. The perturbed trace enters only through its values on \(\Gamma\).
	This proves
	\begin{equation}\label{eq:exact-moment}
		\int_{\partial D} h(\partial_\nu u_0)(\partial_\nu w_g) \dd\sigma=-\int_\Gamma u_\Gamma'[h] g \dd\sigma.
	\end{equation}
	Multiplying both sides of \eqref{eq:exact-moment} by \(\varepsilon\) gives
	\eqref{eq:patch-moment} with
	\(R_\varepsilon(g)=\int_\Gamma
	(  u_\varepsilon-u_0-\varepsilon u_\Gamma'[h])g \dd\sigma\).
	Cauchy--Schwarz and \eqref{eq:local-remainder}, using the continuous
	inclusion \(H^{1/2}(\Gamma)\subset L^2(\Gamma)\), give the asserted bound.
	Equation \eqref{eq:eigen-moment} follows from \eqref{eq:local-shape}
	and \eqref{eq:eigen-remainder}.
	
	Finally, \(\mathcal G_\Gamma=(u_0|_\Gamma)^\perp\) in \(L^2(\Gamma)\).
	The normalization places \(u_\Gamma'[h]|_\Gamma\) in that same subspace.
	Thus vanishing of all the exact adjoint moments implies
	\(u_\Gamma'[h]|_\Gamma=0\), by taking \(g=u_\Gamma'[h]|_\Gamma\)
	in \eqref{eq:exact-moment}. If the eigenvalue moment also vanishes,
	\eqref{eq:local-shape} gives \(\partial_\nu u_\Gamma'[h]=0\) on \(\Gamma\).
	These are the local Cauchy data used in the identifiability argument.
\end{proof}

\begin{proof}[Proof of Theorem \ref{theorem:uniqueness}]
	The first moment in \eqref{eq:h:system} gives \(\lambda'[h]=0\).
	The second and \eqref{eq:exact-moment} show that the trace of
	\(u_\Gamma'[h]\) is orthogonal to \(\mathcal G_\Gamma\). The
	patch normalization in \eqref{eq:local-shape} puts that trace
	in \(\mathcal G_\Gamma\), so it vanishes. The outer Robin
	condition in \eqref{eq:local-shape} then gives zero conormal
	trace on \(\Gamma\). By Lemma~\ref{lemma:UCP}, \(u_\Gamma'[h]=0\)
	in \(A\). Its inner trace therefore gives
	\(h\partial_\nu u_0=0\) on \(\partial D\).
	
	If \(h\) were nonzero, continuity would give a nonempty open inner
	patch on which \(h\) does not vanish. Both the Dirichlet and normal
	traces of \(u_0\) would vanish there. Lemma~\ref{lemma:UCP} would
	then imply \(u_0=0\) in \(A\), contradicting its normalization.
\end{proof}

\subsection{Nonlinear Uniqueness for the Principal Mode}
The nonlinear uniqueness proof uses strict interior positivity of the principal eigenfunction.
\begin{lemma}[Positivity of the principal eigenfunction]\label{lem:positivity:firsteig}
	The principal Steklov--Dirichlet eigenvalue on \(A\) is simple, and its eigenfunction can be chosen strictly positive in \(A\).
\end{lemma}
\begin{proof}
	The argument of \cite[Lemma 1]{HLS} gives simplicity and a
	principal eigenfunction of one sign. Choosing its nonnegative
	representative, the strong minimum principle gives strict
	positivity in the connected domain \(A\).
\end{proof}
\begin{theorem}[Nonlinear identifiability from the principal eigenpair]\label{thm:nonlinear:unique}
	For \(i=1,2\), let \(D_i\) be the inclusion obtained from \(D\) by a sufficiently small \(C^{3,\alpha}(\partial D)\) normal displacement
	\(p_i\), measured along the outward normal \(n_D\) of \(D\). Let \((\lambda_i,u_i)\) be the principal Steklov--Dirichlet eigenpair
	on \(A_i=\Omega\setminus\overline{D_i}\):
	\[
	\Delta u_i=0\quad\text{in } A_i,\qquad
	u_i=0\quad\text{on }\partial D_i,\qquad
	\partial_\nu u_i=\lambda_i u_i\quad\text{on }\partial\Omega
	\]
	normalized by $\langle u_i,u_0\rangle_{L^2(\Gamma)}=\|u_0\|_{L^2(\Gamma)}^2$,
	where \(u_0\) is the reference principal eigenfunction on \(A\).
	If \(\lambda_1=\lambda_2\) and \(u_1|_\Gamma=u_2|_\Gamma\) in this
	normalization, then \(p_1=p_2\) and \(D_1=D_2\).
\end{theorem}

\begin{proof}
	In signed normal coordinates \((y,t)\) in a tubular collar of \(\partial D\) defined by \(x=y+t \mathsf{n}_D(y)\), with \(t>0\)
	on the exterior side, the common exterior \(E=A_1\cap A_2\) is given by \(t>m(y)\), where \(m=\max\{p_1,p_2\}\).
	Choose a smooth cutoff \(\chi(t)\) equal to one near \(t=0\) and zero near the collar ends. For small displacements,
	\[
	(y,t)\longmapsto (y,t+\chi(t)m(y))
	\]
	is strictly increasing on each normal fiber, since
	\(\|\chi'\|_\infty\|m\|_\infty<1\). It maps \(t=0\) to \(t=m(y)\)
	and extends by the identity off the collar to a homeomorphism from \(A\) onto \(E\). Hence \(E\) is connected. The same argument
	with \(p_i\) in place of \(m\) gives connectedness of each \(A_i\).
	
	By Lemma~\ref{lem:positivity:firsteig}, choose the principal modes positive in their respective
	domains. A principal mode cannot have zero trace on an open outer
	patch: the Steklov condition would also give zero normal trace there,
	and Lemma~\ref{lemma:UCP} would force the mode to vanish. Thus the nonnegative
	outer traces of \(u_i\) and \(u_0\) have positive pairing on \(\Gamma\).
	Indeed, the continuous reference trace is positive on a smaller
	open patch, on which the trace of \(u_i\) cannot vanish identically.
	The prescribed scaling therefore preserves positivity.
	
	Set \(z=u_1-u_2\) in \(E\). On \(\Gamma\) the data give
	\(z=0\) and
	\(\partial_\nu z=\lambda_1u_1-\lambda_2u_2=0\).
	Choose a connected component of the outer boundary that meets
	\(\Gamma\) and its thin inner collar, bounded on the inside by an
	inward parallel surface. This collar is contained in \(E\) and has
	\(C^{1,\alpha}\) boundary. Apply Lemma~\ref{lemma:UCP} there on the
	relatively open portion of \(\Gamma\) on its outer boundary. Thus \(z\) vanishes in that collar,
	and interior harmonic unique continuation gives \(z=0\) in \(E\).
	
	If the graph heights differ, interchange the indices if necessary
	so that \(p_1>p_2\) on an open patch of \(\partial D\).
	The corresponding patch of \(\partial D_1\) lies in \(A_2\), and
	its exterior side locally lies in \(E\). Since \(u_1=u_2\) there
	and \(u_1\) has zero Dirichlet trace on that patch, \(u_2\) vanishes
	there by interior continuity. This contradicts positivity in \(A_2\).
	Hence \(p_1=p_2\) and \(D_1=D_2\).
\end{proof}

We briefly review the use of the hypotheses.
The continuation step identifies the eigenfunctions only on the component of the common exterior that meets the observation patch. The collar argument supplies the connectedness needed to reach every differing graph patch. The final contradiction uses the principal mode: a higher eigenfunction may vanish on an interior hypersurface, so equality on the common exterior alone does not complete the same proof for a higher mode.

\subsection{Adjoint Reconstruction and Finite-Dimensional Stability}

To construct compatible tests on \(\Gamma\), choose
\(\psi\in C_c^\infty(\Gamma)\) with
\(\langle\psi,u_0\rangle_{L^2(\Gamma)}\ne0\), which is possible
because \(m_\Gamma>0\). For candidate
\(\varphi_\ell\in C_c^\infty(\Gamma)\), set
\[
g_\ell=\varphi_\ell-
\frac{\langle\varphi_\ell,u_0\rangle_{L^2(\Gamma)}}
{\langle\psi,u_0\rangle_{L^2(\Gamma)}}\,\psi .
\]
Then \(g_\ell\in\mathcal G_\Gamma\). The choice of tests affects
the rank and conditioning of the reconstruction matrix. The
spherical patch experiment uses a projection to remove the same
reference-trace direction.

Let \(\Phi=\operatorname{span}\{\phi_1,\ldots,\phi_M\}
\subset C^{3,\alpha}(\partial D)\), with an
\(L^2(\partial D)\)-orthonormal basis, and choose compatible tests
\(g_1,\ldots,g_L\). Write \(q_g=\partial_\nu w_g|_{\partial D}\)
and let \(\widehat u_\varepsilon\) denote the measured trace
normalized by the patch pairing in Theorem~\ref{thm:simple}.
The moments \eqref{eq:patch-moment} and \eqref{eq:eigen-moment}
give
\[
\begin{aligned}
	K_{0j}&=\int_{\partial D}\phi_j(\partial_\nu u_0)^2\,d\sigma,
	& b_0&=\lambda_\varepsilon-\lambda_0,\\
	K_{\ell j}&=\int_{\partial D}\phi_j(\partial_\nu u_0)q_{g_\ell}\,d\sigma,
	& b_\ell&=-\int_\Gamma(\widehat u_\varepsilon-u_0)g_\ell\,d\sigma
	\quad \forall \ell\in\{1,\dots,L\}.
\end{aligned}
\]
Write \(b^{\rm obs}\) for the vector obtained by using the observed eigenvalue and trace in these formulas. For full-boundary data set \(\Gamma=\partial\Omega\).
In the objective below, \(b\) denotes the signed vector supplied to the solver. The reconstruction \(p_{\rm rec}=\sum_{j=1}^M c_j\phi_j\) is obtained
by minimizing
\begin{equation}\label{eq:tikh}
	\|W^{1/2}(Kc-b)\|_2^2+\alpha_{\mathrm{reg}}\|Rc\|_2^2,\qquad \alpha_{\mathrm{reg}}\ge0,
\end{equation}
where \(K\in\mathbb R^{(L+1)\times M}\), \(W\in\mathbb R^{(L+1)\times(L+1)}\)
is symmetric positive definite, and \(R\in\mathbb R^{q\times M}\)
is a prescribed penalty matrix with \(q\) rows. The normal equations are
\[
(K^TWK+\alpha_{\mathrm{reg}} R^TR)c=K^TWb.
\]
The method estimates \(p\) directly from \(K\) and \(b\).
Smallness of the displacement controls the linearization error.
Quadrature and forward-solver errors enter the data or model error.
If \(p\notin\Phi\), one must also account for the approximation error
and its effect on the data.

\begin{proposition}[Finite-dimensional stability]\label{prop:local}
	Assume the hypotheses of Theorem~\ref{thm:simple}, with the same
	patch normalization and compatible adjoints. Fix a nonzero
	finite-dimensional space \(\Phi\subset C^{3,\alpha}(\partial D)\)
	and an \(L^2(\partial D)\)-orthonormal basis of \(\Phi\).
	\begin{enumerate}
		\item[\textup{(a)}] There exist finitely many tests
		\(g_1,\ldots,g_L\in C_c^\infty(\Gamma)\cap\mathcal G_\Gamma\)
		such that the matrix \(K\) of the eigenvalue moment and the
		selected adjoint moments has full column rank.
		
		\item[\textup{(b)}] Fix a collection of tests as in \textup{(a)} and
		let \(s_\Phi=\sigma_{\min}(K)>0\). Suppose \(p=\varepsilon h\),
		where \(h\) belongs to a fixed bounded subset of \(\Phi\) in the
		\(C^{3,\alpha}\) norm and \(|\varepsilon|\) is sufficiently small.
		If the signed data vector has Euclidean measurement error at most
		\(\eta\), the unweighted ($W=I$), unregularized (\(\alpha_{\mathrm{reg}}=0\)) least-squares reconstruction
		satisfies
		\begin{equation}\label{eq:finitebound}
			\|p_{\rm rec}-p\|_{L^2(\partial D)}
			\le s_\Phi^{-1}\bigl(\eta+C\varepsilon^2\bigr).
		\end{equation}
		The constant \(C\) is uniform on the prescribed bounded set and
		may depend on the selected tests.
	\end{enumerate}
\end{proposition}

\begin{proof}
	\textup{(a)} First, smooth compatible tests are dense in
	\(\mathcal G_\Gamma\) in \(L^2(\Gamma)\). Indeed, for any
	\(g\in\mathcal G_\Gamma\), choose
	\(\varphi_n\in C_c^\infty(\Gamma)\) with \(\varphi_n\to g\) in
	\(L^2(\Gamma)\). Using a fixed bump \(\psi\) with nonzero pairing
	against \(u_0\), as above, gives
	\[
	g_n=\varphi_n-
	\frac{\langle\varphi_n,u_0\rangle_{L^2(\Gamma)}}
	{\langle\psi,u_0\rangle_{L^2(\Gamma)}}\,\psi
	\in C_c^\infty(\Gamma)\cap\mathcal G_\Gamma,
	\qquad g_n\to g\quad\text{in }L^2(\Gamma).
	\]
	The convergence follows because \(g\) is compatible, so the
	correcting coefficient tends to zero. By the adjoint estimate in
	Theorem~\ref{thm:adjoint}, each adjoint moment depends continuously
	on the test in \(L^2(\Gamma)\). Thus if a displacement in \(\Phi\)
	annihilates the eigenvalue moment and every smooth compatible
	moment, it also annihilates every compatible moment.
	Theorem~\ref{theorem:uniqueness} then forces the displacement to be
	zero. Starting with the kernel of the eigenvalue row, choose a
	smooth compatible test that is nonzero on the current kernel
	whenever that kernel is nontrivial. Each such test lowers its
	dimension by one. The process terminates after finitely many
	tests and yields full column rank.
	
	\textup{(b)} Let \(c\) be the coefficient vector of \(p\), and let
	\(b^{\rm obs}\) be the observed signed data vector. The displacement
	identities in Theorem~\ref{thm:adjoint}, applied to the fixed tests,
	give
	\[
	b^{\rm obs}=Kc+r+e,\qquad
	\|r\|_2\le C\varepsilon^2,\qquad \|e\|_2\le\eta .
	\]
	Here \(r\) is the signed linearization remainder, uniform for the
	bounded perturbation class, and \(e\) is the difference between the
	observed signed vector and the signed vector formed from the noise-free
	finite-amplitude eigenvalue and trace.
	Since \(K\) has full column rank, \(c_{\rm rec}=K^\dagger b^{\rm obs}\)
	and \(K^\dagger K=I\). Orthonormality of the basis and
	\(\|K^\dagger\|_{2\to2}=s_\Phi^{-1}\) give
	\[
	\|p_{\rm rec}-p\|_{L^2(\partial D)}
	=\|c_{\rm rec}-c\|_2
	=\|K^\dagger(r+e)\|_2
	\le s_\Phi^{-1}\bigl(\eta+C\varepsilon^2\bigr).
	\]
\end{proof}

The stability constant depends on the reference geometry, the patch,
the tests, and the trial space. The full-boundary spherical estimate
in Section~\ref{sec:sphere} concerns a different data norm and gives
an infinite-dimensional rate under a Sobolev prior.

\section{Shape Sensitivity and Linearized Identifiability for Neumann Inclusions}\label{sec:neumann}

Replace the inner Dirichlet condition in \eqref{eq:evp} by
\(\partial_\nu u_\varepsilon^N=0\) on \(\partial D_\varepsilon\).
The pullback uses the same forms \(a_\varepsilon,b\), now on
\(H^1(A)\). The trace Poincar\'e inequality
\[
\|z\|_{H^1(A)}^2\le C\bigl(\|\nabla z\|_{L^2(A)}^2+
\|z\|_{L^2(\partial\Omega)}^2\bigr)
\]
makes \(a_\varepsilon+b\) uniformly coercive for small \(|\varepsilon|\).
The compact operator \((a_\varepsilon+b)^{-1}b\) has nonzero
eigenvalues \((1+\lambda_\varepsilon^N)^{-1}\).
Testing with \(1\) gives zero outer mean for positive modes.
Connectedness makes the zero eigenspace one-dimensional.
Its normalized constant eigenfunction \(|\partial\Omega|^{-1/2}\)
gives the same outer data for every admissible inner shape.

\subsection{Neumann Regularity and Shape Derivatives}

Let \((\lambda_0^N,u_0^N)\) be a positive simple eigenpair, normalized
by \(\|u_0^N\|_{L^2(\partial\Omega)}=1\).
On any nonempty relatively open \(\Gamma\subset\partial\Omega\), set
\[
\langle u_\varepsilon^N,u_0^N\rangle_{L^2(\Gamma)}
=\|u_0^N\|_{L^2(\Gamma)}^2 .
\]
The right-hand side is positive: a zero trace on \(\Gamma\) gives a
zero normal trace by the Steklov condition, and boundary unique
continuation then contradicts normalization. The bordered argument
of Lemma~\ref{lem:analytic} applies with \(a_0+b\) as the coercive
form. It gives an analytic branch in \(\mathbb R\times H^1(A)\)
and the quadratic remainders
\eqref{eq:eigen-remainder}--\eqref{eq:local-remainder}
with the corresponding Neumann quantities.

\begin{lemma}[Regularity and weak shape derivative]
	Write \(u=u_0^N\), \(\lambda=\lambda_0^N\),
	\(U_\varepsilon=u_\varepsilon^N\circ F_\varepsilon\), and
	\(v=\dot U-V\cdot\nabla u\), with dots denoting derivatives at zero.
	Then \(u\in H^3(A)\), \(v\in H^2(A)\), and
	\begin{equation}\label{eq:SNderiv}
		\begin{cases}
			\Delta v=0 & \text{in }A,\\
			(\partial_\nu-\lambda)v=\dot\lambda u & \text{on }\partial\Omega,\\
			\partial_\nu v=-\operatorname{div}_\tau(h\nabla_\tau u)
			& \text{on }\partial D,\\
			\langle v,u\rangle_{L^2(\Gamma)}=0 .
		\end{cases}
	\end{equation}
	Here \(\nabla_\tau\) and \(\operatorname{div}_\tau\) are the surface
	gradient and divergence on \(\partial D\), and \(\nu\) points out of \(A\).
	Moreover,
	\begin{equation}\label{eq:SNlambda}
		(\lambda^N)'[h]=\dot\lambda
		=-\int_{\partial D} h|\nabla_\tau u|^2\dd\sigma .
	\end{equation}
\end{lemma}
\begin{proof}
	The outer Neumann datum \(\lambda u\) has \(H^{1/2}\) regularity
	by the \(H^1\) trace theorem. The inner datum is zero.
	The \(H^2\) Neumann estimate improves the outer datum to \(H^{3/2}\),
	and a second application at the next order gives \(u\in H^3(A)\)
	\cite[Ch.~4]{McLean}. In particular,
	\(u|_{\partial D}\in H^{5/2}(\partial D)\).
	
	Differentiate the pulled-back weak equation. Since \(V\) vanishes
	near \(\partial\Omega\), for \(\phi\in H^1(A)\) this gives
	\[
	\int_A(\nabla v+J)\cdot\nabla\phi\dd x
	=\int_{\partial\Omega}(\lambda v+\dot\lambda u)\phi\dd\sigma,
	\qquad J=C'_0\nabla u+\nabla(V\cdot\nabla u),
	\]
	where \(C'_0=(\operatorname{div}V)I-DV-DV^T\) is the derivative
	of the pullback coefficient, where \((DV)_{ij}=\partial_j V_i\). With summation over \(j\),
	\[
	J_i=\partial_j\bigl((\partial_i u)V_j-V_i\partial_j u\bigr),
	\]
	using \(\Delta u=0\). For smooth tests, integration by parts leaves
	only a boundary term because the tensor in parentheses is
	antisymmetric. The outer term vanishes, and on \(\partial D\)
	one has \(\partial_\nu u=0\) and \(V\cdot\nu=-h\). Thus
	\[
	\int_A J\cdot\nabla\phi\dd x
	=-\int_{\partial D} h\nabla_\tau u\cdot\nabla_\tau\phi\dd\sigma
	=\int_{\partial D}\operatorname{div}_\tau(h\nabla_\tau u)\phi\dd\sigma .
	\]
	The surface divergence belongs to \(H^{1/2}(\partial D)\) by the
	trace regularity of \(u\), so the identity extends by density.
	The weak equation now gives the conormal conditions in
	\eqref{eq:SNderiv}. The outer datum \(\lambda v+\dot\lambda u\)
	belongs to \(H^{1/2}(\partial\Omega)\), so the \(H^2\) Neumann
	estimate gives \(v\in H^2(A)\). The last condition in
	\eqref{eq:SNderiv} follows by differentiating the patch pairing.
	Finally, test the equation for \(v\) with \(u\) and the reference
	equation with \(v\). Subtract and use \(\|u\|_{L^2(\partial\Omega)}=1\)
	to obtain \eqref{eq:SNlambda}.
\end{proof}

\subsection{Local Adjoint Measurements and Linearized Identifiability}

Let \(g\in L^2(\Gamma)\) be orthogonal to \(u_0^N|_\Gamma\) and
extend it by zero to \(\partial\Omega\). For a positive simple mode,
the Fredholm alternative gives a unique harmonic adjoint \(w_g^N\)
with \(\partial_\nu w_g^N=0\) on \(\partial D\),
\((\partial_\nu-\lambda_0^N)w_g^N=g\) on \(\partial\Omega\), and
\(\ip{w_g^N}{u_0^N}=0\). Homogeneous Neumann regularity in an inner
collar gives its tangential trace gradient there. Testing with \(1\)
gives \(-\lambda_0^N\int_{\partial\Omega}w_g^N\dd\sigma=\int_\Gamma g\dd\sigma\).
Green's identity with \(v\), followed by the trace remainder, gives
for \(p=\varepsilon h\)
\begin{equation}\label{eq:SNdual}
	\int_{\partial D}p\,\nabla_\tau u_0^N\cdot\nabla_\tau w_g^N\dd\sigma
	=\int_\Gamma(u_\varepsilon^N-u_0^N)g\dd\sigma
	+O(\varepsilon^2)\|g\|_{L^2(\Gamma)}.
\end{equation}
The eigenvalue shift gives the moment in \eqref{eq:SNlambda} with
\(p\) in place of \(h\), up to \(O(\varepsilon^2)\).

\begin{theorem}[Linearized Steklov--Neumann identifiability]\label{thm:injectivity}
	For any positive simple Steklov--Neumann eigenvalue and any nonempty
	relatively open \(\Gamma\subset\partial\Omega\), the map
	\(h\mapsto ((\lambda^N)'[h],v|_\Gamma)\) is injective on real
	\(C^{3,\alpha}(\partial D)\).
\end{theorem}
\begin{proof}
	If both data vanish, \eqref{eq:SNderiv} gives zero Dirichlet and
	conormal traces on \(\Gamma\). Lemma~\ref{lemma:UCP} gives
	\(v=0\) in \(A\).
	With \(f=u_0^N|_{\partial D}\in H^{5/2}(\partial D)\), the inner
	condition therefore gives \(\operatorname{div}_\tau(h\nabla_\tau f)=0\).
	
	The Sobolev chain rule gives
	\(\nabla_\tau |h|=\operatorname{sgn}(h)\nabla_\tau h\) almost
	everywhere, with \(\operatorname{sgn}(0)=0\). The product rule gives
	\[
	\operatorname{div}_\tau(|h|\nabla_\tau f)
	=\operatorname{sgn}(h)\operatorname{div}_\tau(h\nabla_\tau f)=0 .
	\]
	Testing on the closed inner boundary with \(f\) yields
	\(\int_{\partial D}|h|\,|\nabla_\tau f|^2\dd\sigma=0\).
	If \(h\) is nonzero, it is bounded away from zero in absolute value
	on a connected open inner patch. There \(\nabla_\tau f=0\) almost
	everywhere, so \(f\) is constant as a Sobolev trace. Subtract this
	constant from \(u_0^N\) and apply Lemma \ref{lemma:UCP} using the homogeneous inner conormal condition. The resulting
	constant eigenfunction must be zero because \(\lambda_0^N>0\),
	contradicting normalization.
\end{proof}

On a fixed finite-dimensional smooth trial space, the eigenvalue
moment and compatible trace tests separate displacements because
\(v|_\Gamma\) is itself compatible. A finite full-rank collection
can therefore be selected as in Proposition~\ref{prop:local}.
For unweighted least squares its smallest singular value controls
the linearized noise amplification in an \(L^2(\partial D)\)-orthonormal
basis.

\section{Spherical Harmonics and Full-Boundary Conditional Stability}\label{sec:sphere}
In this section, we establish the conditional stability with respect to full-boundary data of a principal eigenpair for the case where the reference domain is the concentric shell in 3D. We set $d=3$, $\Om=B_1(0)$, and $D=B_a(0)$ for some $0<a<1$.
The principal eigenpair has the explicit form below~\cite[Remark 2.1]{VS:2020:EDC}.
\begin{equation}\label{eq:radial}
	\lambda_0=\frac{a}{1-a}\quad\mbox{and}\quad
	u_0(x)=\frac{a}{\sqrt{4\pi}(1-a)}\left(\frac{1}{a}-\frac{1}{|x|}\right)\quad \forall x\in \Om\backslash\overline{D},\quad\mbox{thus}\quad\partial_\nu u_0=-\frac1{\sqrt{4\pi}a(1-a)}.
\end{equation}
For $\mathbb{S}^2$, the appropriate analogue of a Fourier basis is a real
orthonormal spherical-harmonic basis $Y_{\ell m}$,
$-\ell\leq m\leq\ell$, with $Y_{00}=1/\sqrt{4\pi}$.
We expand the displacement pulled back to the unit sphere as
$$\widetilde{p}(\theta):=p(a\theta)=\sum_{\ell=0}^\infty\sum_{m=-\ell}^\ell p_{\ell m}Y_{\ell m}(\theta)\quad\forall \theta\in\mathbb{S}^2.$$
These coefficients use $\dd\theta$, not the area measure $a^2\dd\theta$ on $\partial B_a$. Consequently $\norm p_{L^2(\partial B_a)}=a\norm{\widetilde{p}}_{L^2(\mathbb{S}^2)}$.
For a function on $\mathbb{S}^2$ we use the spectral Sobolev norm
\[
\norm{\widetilde{p}}_{H^s(\mathbb{S}^2)}^2=
\sum_{\ell=0}^\infty\sum_{m=-\ell}^{\ell}
(1+\ell(\ell+1))^s|p_{\ell m}|^2.
\]
The coefficient map $\widetilde{p}\mapsto \{p_{\ell m}\}_{\ell,m}$ is unitary from $L^2(\mathbb{S}^2)$ onto
$\ell^2(\{(\ell,m)\,:\,\ell\geq0,-\ell\leq m\leq\ell\})$.

For $\ell\geq0$ and $a\le r\le 1$, define
\[
R_\ell(r)=\frac{r^\ell-a^{2\ell+1}r^{-\ell-1}}{1-a^{2\ell+1}}
\quad\mbox{and}\quad
\sigma_\ell=R_\ell'(1)
=\frac{\ell+(\ell+1)a^{2\ell+1}}{1-a^{2\ell+1}}.
\]
Then $R_\ell(a)=0$ and $R_\ell(1)=1$. In particular, $\sigma_0=\lambda_0$. For all $\ell\geq1$ and $|m|\le\ell$, $g=Y_{\ell m}$ is compatible on the full boundary because $u_0|_{\mathbb{S}^2}=Y_{00}$ and $\langle Y_{\ell m},Y_{00}\rangle_{L^2(\mathbb{S}^2)}=0$,
and the unique solution of the adjoint equation \eqref{eq:adjoint} is
\[w_g(r,\theta):=
\frac{R_\ell(r)Y_{\ell m}(\theta)}{\sigma_\ell-\lambda_0}
\quad\forall(r,\theta)\in(a,1)\times\mathbb{S}^2\]
in spherical coordinates.
The denominator is positive since the principal eigenvalue is simple and strictly smaller than the angular modes.
We therefore have, for $\ell\geq1$ and $g=Y_{\ell m}$,
\begin{equation}\label{eq:adj:explicit}
	\p_\nu w_g(\theta)=-\frac{(2\ell+1)a^{\ell-1}}{(1-a^{2\ell+1})(\sigma_\ell-\lambda_0)}Y_{\ell m}(\theta)\quad\forall \theta\in\mathbb{S}^2.
\end{equation}
Substituting these expressions into the moment identities of Theorem~\ref{thm:adjoint} gives a diagonal reconstruction formula.

\begin{proposition}[Diagonal reconstruction]\label{prop:diag}
	Set
	\begin{align}
		\kappa_0&=\frac1{\sqrt{4\pi}(1-a)^2},\label{eq:k0}\\
		\kappa_\ell&=\frac{(2\ell+1)a^\ell}
		{\sqrt{4\pi}(1-a)(1-a^{2\ell+1})(\sigma_\ell-\lambda_0)},
		\quad\forall \ell\geq1.\label{eq:kl}
	\end{align}
	Let $t_{\ell m}=\langle u_\varepsilon-u_0,
	Y_{\ell m}\rangle_{L^2(\mathbb{S}^2)}$. With full-boundary normalization, \(\Gamma=\partial\Omega\), we obtain
	\begin{equation}\label{eq:spectraldata}
		\lambda_\varepsilon-\lambda_0=\kappa_0p_{00}+O(\varepsilon^2)
		\quad\mbox{and}\quad
		t_{\ell m}=-\kappa_\ell p_{\ell m}+O(\varepsilon^2)
		\quad\forall\ell\geq1\ \forall m\in\{-\ell,-\ell+1,\dots,\ell\}.
	\end{equation}
	The constant outer trace coefficient vanishes identically. Thus the trace coefficients \(t_{\ell m}\) alone cannot recover \(p_{00}\).
\end{proposition}
\begin{proof}
	We apply Theorem \ref{thm:adjoint}.
	The eigenvalue formula \eqref{eq:eigen-moment} and the expression \eqref{eq:radial} give
	$$\lambda_\varepsilon-\lambda_0+O(\varepsilon^2) = \int_{\partial B_a(0)}p(\partial_\nu u_0)^2\dd\sigma = a^2(\partial_\nu u_0)^2\int_{\mathbb{S}^2}\widetilde{p}(\theta)\dd\theta=\kappa_0p_{00}.$$
	For \(g=Y_{\ell m}\) with \(\ell\geq1\), the trace moment
	\eqref{eq:patch-moment}, the formulas \eqref{eq:radial} and
	\eqref{eq:adj:explicit}, and spherical-harmonic orthogonality give
	$$-t_{\ell m}+O(\varepsilon^2)=\int_{\partial B_a(0)}p (\partial_\nu u_0)(\partial_\nu w_g)\dd\sigma=\kappa_\ell \int_{\mathbb{S}^2} \widetilde{p}(\theta)Y_{\ell m}(\theta)\dd\theta=\kappa_\ell p_{\ell m}.$$
	Finally, $\ip{u_\varepsilon-u_0}{u_0}=0$ and
	$u_0|_{\mathbb{S}^2}=Y_{00}$ make the constant trace coefficient zero.
\end{proof}

For fixed $a$, $\sigma_\ell=\ell+o(1)$ and consequently
\begin{equation}\label{eq:decay}
	\kappa_\ell\sim\frac{2a^\ell}{\sqrt{4\pi}(1-a)}
	\quad\text{as }\ell\to\infty.
\end{equation}
In particular there are positive constants $c_a,C_a$ such that
$c_a a^\ell\leq\kappa_\ell\leq C_a a^\ell$ for all $\ell\geq1$.
The inverse amplifies high-degree noise exponentially.

The diagonal coefficient map below is the bounded compact extension
of the smooth-displacement linearization from $L^2(\mathbb{S}^2)$ into
$\R\oplus\ell^2(\{(\ell,m):\ell\geq1\})$, with the standard
product norm. The first coordinate is the eigenvalue perturbation;
the other coordinates are the sign-adjusted outer trace coefficients.
This algebraic extension permits a Sobolev prior with any $s>0$ in
the linearized estimate, without asserting that all such functions
are admissible nonlinear boundary graphs.

\begin{theorem}[Full-boundary cutoff stability]\label{thm:stability}
	For the diagonal linearized map with data
	$d=(\kappa_0p_{00},\{\kappa_\ell p_{\ell m}\}_{\ell\geq1})$,
	assume $\norm{\widetilde{p}}_{H^s(\mathbb{S}^2)}\leq M$, $s>0$, and
	\(\|d^\eta-d\|_{\mathbb R\oplus\ell^2}\leq\eta\), where the scalar
	coordinate corresponds to the eigenvalue shift. The degree-$L$ cutoff inverse
	satisfies, for $L\geq1$,
	\begin{equation}\label{eq:cutoff}
		\norm{p_L^\eta-\widetilde{p}}_{L^2(\mathbb{S}^2)}
		\leq C_a a^{-L}\eta+C_sM(L+1)^{-s}.
	\end{equation}
	For $\eta/M$ sufficiently small, choosing
	$L=\lfloor\log(M/\eta)/(2|\log a|)\rfloor$ yields
	\[
	\norm{p_L^\eta-\widetilde{p}}_{L^2(\mathbb{S}^2)}
	\leq C_{a,s}M\,[\log(M/\eta)]^{-s}.
	\]
	There is no bounded inverse in the unweighted $L^2$ data and
	displacement norms.
\end{theorem}
\begin{proof}
	Write \(d^\eta_{\ell m}\) for the observed data coefficients, with
	\(d^\eta_{00}\) denoting the eigenvalue coordinate. The cutoff
	reconstruction has coefficients \(d^\eta_{\ell m}/\kappa_\ell\)
	for \(\ell\leq L\) and zero coefficients for \(\ell>L\).
	Since \(d_{\ell m}=\kappa_\ell p_{\ell m}\), Parseval's identity gives
	\[
	\|p_L^\eta-\widetilde p\|_{L^2(\mathbb S^2)}^2
	=\sum_{\ell=0}^{L}\sum_{m=-\ell}^{\ell}
	\frac{|d^\eta_{\ell m}-d_{\ell m}|^2}{\kappa_\ell^2}
	+\sum_{\ell>L}\sum_{m=-\ell}^{\ell}|p_{\ell m}|^2 .
	\]
	The two terms represent propagated noise and the discarded Sobolev
	tail, respectively. The two-sided bounds following
	\eqref{eq:decay}, together with the fixed factor \(\kappa_0>0\),
	control the first term. The Sobolev prior controls the second:
	\[
	\begin{aligned}
		\sum_{\ell=0}^{L}\sum_{m=-\ell}^{\ell}
		\frac{|d^\eta_{\ell m}-d_{\ell m}|^2}{\kappa_\ell^2}
		&\leq C_a^2a^{-2L}\eta^2,\\
		\sum_{\ell>L}\sum_{m=-\ell}^{\ell}|p_{\ell m}|^2
		&\leq M^2\bigl(1+(L+1)(L+2)\bigr)^{-s}
		\leq M^2(L+1)^{-2s}.
	\end{aligned}
	\]
	Taking square roots gives \eqref{eq:cutoff}. For the prescribed
	choice of \(L\), put \(T=\log(M/\eta)\). Since
	\(L\leq T/(2|\log a|)<L+1\), we obtain
	\[
	a^{-L}\eta\leq\sqrt{M\eta}=M e^{-T/2},
	\qquad
	M(L+1)^{-s}\leq M(2|\log a|)^sT^{-s}.
	\]
	For large \(T\), \(e^{-T/2}\leq C_sT^{-s}\), which gives the
	logarithmic estimate. Finally, a unit-norm spherical harmonic of
	degree \(\ell\) has data norm \(\kappa_\ell\to0\) by
	\eqref{eq:decay}. A bounded inverse would bound its unit
	displacement norm by a fixed multiple of that vanishing data norm.
\end{proof}

For actual nonlinear measurements, \eqref{eq:eigen-remainder} and
\eqref{eq:local-remainder} with \(\Gamma=\partial\Omega\), together
with Parseval's identity, bound the full data remainder by
\(C\varepsilon^2\) in \(\mathbb R\oplus\ell^2\). Thus \eqref{eq:cutoff} holds with an effective error
$\eta+C\varepsilon^2$ in place of $\eta$, provided $h$ remains in
a bounded smooth class. A larger cutoff eventually amplifies the
model error as well as the measurement noise.

\paragraph{Diagonal regularization.}
We use a Tikhonov penalty~\cite[Chapter 8]{EHN} with weights
\(r_\ell^2=(1+\ell(\ell+1))^s\). Since \(\ell(\ell+1)\) is the
eigenvalue of \(-\Delta_{\mathbb S^2}\) on degree-\(\ell\) harmonics,
these weights give the squared \(H^s(\mathbb S^2)\) norm of the
displacement, consistent with the prior in Theorem~\ref{thm:stability}.
The full-boundary regularized coefficients are
\[
p_{00}^{\alpha_{\mathrm{reg}}}=\frac{\kappa_0d_{00}}{\kappa_0^2+\alpha_{\mathrm{reg}} r_0^2},
\qquad
p_{\ell m}^{\alpha_{\mathrm{reg}}}=
\frac{\kappa_\ell d_{\ell m}}{\kappa_\ell^2+\alpha_{\mathrm{reg}} r_\ell^2},
\qquad d_{\ell m}=-t_{\ell m}\ (\ell\geq1).
\]
The measured eigenvalue supplies $d_{00}=\lambda_\varepsilon-\lambda_0$.
The penalty suppresses modes for which \(\kappa_\ell^2\ll\alpha_{\mathrm{reg}} r_\ell^2\), where direct inversion would amplify noise by \(1/\kappa_\ell\). Section~\ref{sec:localnumerics} uses the same weights with \(s=2\) on the nonconstant modes in the local fit.

\section{Numerical Examples}
\label{sec:numerics}

We test the linearized Dirichlet reconstruction on data from the principal eigenpair of the perturbed shell. The forward solver uses the conforming Rayleigh--Ritz approximation for the exact perturbed geometry. The full-boundary experiment tests the reconstruction error as the perturbation amplitude decreases.
The local-patch experiment tests conditioning and noise sensitivity as the observed cap shrinks.

\subsection{Variational Forward Solver}

For simplicity, throughout the numerical tests, we only consider the perturbed shell invariant under rotations about its axis.
Each such rotation preserves the positive, normalized principal
eigenfunction by simplicity. Thus the principal eigenfunction is axisymmetric, and the minimum of the Rayleigh quotient is attained in the axisymmetric energy space.

Let \(x=\cos\theta\), write the inner boundary as \(r=\rho(x)\), and
set \(g(x)=1-\rho(x)>0\). Write \(r=\rho(x)+\xi g(x)\), so
\(\xi=0\) and \(\xi=1\) correspond to the inner and outer boundaries.
For a pulled-back axisymmetric function \(U(\xi,x)\), the energy and
squared outer trace norm are
\[
\begin{aligned}
	a(U,U)&=2\pi\int_{-1}^{1}\int_0^1
	\left[
	\frac{r^2}{g}|U_\xi|^2+
	g(1-x^2)\left|U_x-\frac{(1-\xi)\rho'}{g}U_\xi\right|^2
	\right]\,d\xi\,dx,\\
	b(U,U)&=2\pi\int_{-1}^{1}|U(1,x)|^2\,dx .
\end{aligned}
\]
We minimize \(a(U,U)/b(U,U)\) over the tensor-product space with
radial basis \(\xi\) and \(\xi(1-\xi)P_j(2\xi-1)\), \(0\le j<m\),
and angular basis \(Y_{\ell,0}(x)=\sqrt{(2\ell+1)/(4\pi)}P_\ell(x)\),
\(0\le\ell\le m\), where \(P_j\) denotes the Legendre polynomial of
degree \(j\). Every trial function vanishes on the inner boundary.
Eliminating the coefficients of the radial functions with zero outer trace gives a
symmetric Schur complement in the outer harmonic coefficients. We
take its lowest eigenpair and rescale the outer trace so that its
\(Y_{0,0}\) coefficient is one.

For a non-axisymmetric radial graph \(r=\rho(\omega)\), \(\omega\in\mathbb S^2\),
one can use the same radial basis with all spherical harmonics \(Y_{\ell,k}\),
\(0\le\ell\le m\), \(-\ell\le k\le\ell\). Assemble the pulled-back
Dirichlet energy with the full spherical gradient and integrate over
\(\mathbb S^2\). The same Schur complement gives the Ritz eigenproblem.

The trial spaces are dense in the axisymmetric energy space, so
standard Rayleigh--Ritz theory gives convergence of the principal
eigenvalue and its normalized outer trace as \(m\to\infty\)
(cf.~\cite{AT:2004:IFESEP} for Steklov finite-element convergence).
We use \(2m+12\) Gauss--Legendre nodes per coordinate. For the
polynomial graphs below, the outer trace norm is integrated exactly
and the energy quadrature error tends to zero relative to the exact
energy, uniformly on the trial space, by polynomial approximation of
\(1/g\) and Gauss exactness. Thus the computed eigenvalue and trace
also converge.

\subsection{Full-Boundary Reconstruction}

Take \(a=1/2\) and
\[
\rho(x)=a+\varepsilon h(x),
\qquad h=0.30Y_{00}+0.50Y_{20}-0.25Y_{30}.
\]
At \(\varepsilon=0.04\), Table~\ref{tab:ritz-refinement} reports
the differences from the \(m=16\) eigenvalue and outer trace.
Increasing the quadrature at \(m=16\) from 44 to 56 points per
coordinate changed the trace by \(6.5\times10^{-15}\) in
\(L^2(\mathbb{S}^2)\). The eigenvalue was unchanged at the precision of
the calculation.

\begin{table}[htbp]
	\centering
	\begin{tabular}{@{}cccc@{}}
		\toprule
		\(m\)&Eigenvalue&Eigenvalue difference&Trace difference\\
		\midrule
		4 &1.014178605879&\(1.03\times10^{-6}\)&\(2.15\times10^{-5}\)\\
		6 &1.014177581473&\(1.54\times10^{-9}\)&\(1.99\times10^{-7}\)\\
		8 &1.014177579934&\(2.19\times10^{-12}\)&\(8.98\times10^{-9}\)\\
		10&1.014177579932&\(3.33\times10^{-15}\)&\(2.23\times10^{-10}\)\\
		12&1.014177579932&\(8.88\times10^{-16}\)&\(4.17\times10^{-12}\)\\
		16&1.014177579932&\multicolumn{2}{c}{reference}\\
		\bottomrule
	\end{tabular}
	\caption{Forward refinement at \(\varepsilon=0.04\).
		Differences use the \(m=16\) solution as a numerical reference.
		Trace differences are measured in \(L^2(\mathbb{S}^2)\).}
	\label{tab:ritz-refinement}
\end{table}

Using the \(m=16\) forward data, we reconstruct the displacement
coefficients through degree six by the diagonal formulas of
Section~5. We use the eigenvalue shift for the constant coefficient
and the outer trace for the remaining coefficients, without noise
or regularization. Table~\ref{tab:conv} reports the absolute
\(L^2(\mathbb{S}^2)\) coefficient error and its observed order under halving
\(\varepsilon\). The order approaches two, in agreement with the
quadratic linearization error at a fixed cutoff.

\begin{table}[htbp]
	\centering
	\begin{tabular}{@{}cccc@{}}
		\toprule
		\(\varepsilon\)&Eigenvalue&Absolute error&Order\\
		\midrule
		0.0400&1.014177579932&\(1.423918\times10^{-3}\)&--\\
		0.0200&1.006930232427&\(3.534174\times10^{-4}\)&2.01042\\
		0.0100&1.003425198347&\(8.801057\times10^{-5}\)&2.00562\\
		0.0050&1.001702591899&\(2.195827\times10^{-5}\)&2.00291\\
		0.0025&1.000848791092&\(5.483936\times10^{-6}\)&2.00148\\
		\bottomrule
	\end{tabular}
	\caption{Full-boundary reconstruction through degree six using
		the perturbed forward eigenpair.}
	\label{tab:conv}
\end{table}

\paragraph{Constant displacement check.}
For a constant displacement \(c\), the exact principal eigenvalue
is \((a+c)/(1-a-c)\), while the normalized outer trace stays
constant. At \(a=0.5\) and \(c=0.002\), the computed eigenvalue
is \(1.008032128514056\), within \(4.4\times10^{-16}\) of
\(0.502/0.498\). The nonconstant trace coefficients have Euclidean
norm \(5.1\times10^{-15}\).

\subsection{Local-Patch Conditioning and Noise Sensitivity}
\label{sec:localnumerics}

We use the same nonspherical perturbation with \(\varepsilon=0.01\)
and retain the outer trace on the cap \(x>x_*\).
The inverse fits degrees one through three from the trace and uses
the eigenvalue shift for the mean. On each cap we use 120
Gauss--Legendre nodes and normalize the quadrature weights after
multiplication by
\[
\chi(x)=\exp\!\left(-\frac{(1-x_*)^2}{(x-x_*)(1-x)}\right)
\qquad (x_*<x<1).
\]
Let \(x_i\) be the cap nodes and \(u_i=u_\varepsilon(1,x_i)\).
Write \(W_i\) for the normalized weights and define
\(\|y\|_W^2=\sum_iW_i|y_i|^2\) for a sampled vector \(y\).
Set \(z_i=\sqrt{W_i}\), \(Q=I-zz^T\), and \(u_*=(4\pi)^{-1/2}\).
For the sampled trace \(u\), put
\[
\widehat u_i=\frac{u_*u_i}{\sum_j W_j u_j},
\qquad d_{\mathrm{cap}}=Q\operatorname{diag}(z)(\widehat u-u_*\boldsymbol 1),
\qquad B_{\cdot\ell}=-\kappa_\ell Q\operatorname{diag}(z)Y_{\ell0},
\quad 1\le\ell\le3 ,
\]
where each \(Y_{\ell0}\) in the matrix is evaluated at the cap
nodes. This projects away the reference-trace direction and gives
the first-order model \(d_{\mathrm{cap}}\approx B(p_{10},p_{20},p_{30})^T\).
The mean coefficient is estimated by \((\lambda_\varepsilon-\lambda_0)/\kappa_0\).

Before normalizing the trace, we add noise proportional to
\(\sin(17\theta)+0.4\cos(29\theta)\) with weighted norm
\(10^{-6}\) times the weighted norm of the sampled trace.
The eigenvalue is kept noiseless. We compare least squares with
Tikhonov regularization using \(\alpha_{\mathrm{reg}}=10^{-10}\) and
penalty weights \(r_\ell^2=(1+\ell(\ell+1))^2\).

\begin{table}[htbp]
	\centering
	\small
	\begin{tabular}{@{}cccccc@{}}
		\toprule
		Cap fraction & \(s_{\min}\) & Noiseless & Noisy & Regularized
		& Noise shift\\
		\midrule
		1.000&\(1.460{\times}10^{-2}\)&\(5.764{\times}10^{-5}\)&\(5.766{\times}10^{-5}\)&\(5.758{\times}10^{-5}\)&\(2.109{\times}10^{-8}\)\\
		0.875&\(9.336{\times}10^{-3}\)&\(5.502{\times}10^{-5}\)&\(5.464{\times}10^{-5}\)&\(5.440{\times}10^{-5}\)&\(8.468{\times}10^{-7}\)\\
		0.625&\(2.622{\times}10^{-3}\)&\(6.374{\times}10^{-5}\)&\(6.103{\times}10^{-5}\)&\(5.639{\times}10^{-5}\)&\(3.899{\times}10^{-6}\)\\
		0.375&\(4.139{\times}10^{-4}\)&\(9.544{\times}10^{-5}\)&\(2.496{\times}10^{-4}\)&\(4.660{\times}10^{-4}\)&\(3.349{\times}10^{-4}\)\\
		\bottomrule
	\end{tabular}
	\caption{Local reconstruction in the normalized weighted observation
		norm. The cap fraction is \((1-x_*)/2\), and \(s_{\min}\) is the smallest singular value of \(B\). The three reconstruction errors are measured against the true displacement in \(L^2(\mathbb{S}^2)\).
		The regularized estimate uses the noisy data. The noise shift is the
		distance between the noisy and noiseless unregularized estimates.}
	\label{tab:local}
\end{table}

As the cap shrinks, the smallest singular value decreases and
the noise shift increases for the prescribed noise profile.
On some intermediate caps, the total error decreases after noise is added, indicating partial cancellation with the linearization error. On the smallest cap, the fixed penalty increases the error. These comparisons concern the stated three-mode trial space and the normalized, patch-dependent observation weights.

\section{Conclusions and Open Questions}\label{sec:scope}
For a known Dirichlet inclusion, one simple eigenvalue and its
locally normalized eigenfunction trace determine the displacement
at first order. On a fixed finite-dimensional trial space, the
smallest singular value of the measurement matrix controls
reconstruction stability. For the principal Dirichlet mode, the
same data give nonlinear uniqueness among sufficiently small
normal graphs. For a Neumann inclusion, every positive simple
mode gives linearized injectivity.

In a three-dimensional concentric shell, the linearized full-boundary
map gives logarithmic conditional stability under a Sobolev prior.
The numerical experiments illustrate the quadratic reconstruction
error at a fixed cutoff and the sensitivity to noise on smaller
observation patches. Quantitative infinite-dimensional stability
from local data and recovery from multiple eigenspaces remain
open questions in this setting.

\section*{Acknowledgement}
The author thanks Professor Bangti Jin of The Chinese University of Hong Kong for his support and advice throughout this project.
The author also thanks Professor Habib Ammari for hosting a short visit to ETH Zurich and for valuable discussions during the project.

\bibliographystyle{abbrv}
\bibliography{reference}

\end{document}